\documentclass[12pt]{amsart}

\usepackage{a4,amsmath,amssymb,amsfonts}
\usepackage[all]{xy}

\newcounter{lemma}
\newtheorem{localtheorem}[lemma]{Theorem}

\begin{document}
\baselineskip=15.5pt

\renewcommand{\phi}{\varphi}

\newcommand{\dual}{^\lor}
\newcommand{\Ext}{{\rm Ext}}
\newcommand{\Ecal}{{\mathcal E}}
\newcommand{\Hom}{{\rm Hom}}
\newcommand{\eps}{\varepsilon}
\newcommand{\inv}{^{-1}}
\newcommand{\Ocal}{{\mathcal O}}
\newcommand{\rarpa}[1]{\stackrel{#1}{\longrightarrow}}
\newcommand{\rk}{{\rm rk}}
\newcommand{\Spec}{{\rm Spec}}
\newcommand{\SU}{{\rm SU}}
\newcommand{\pdop}{{\mathbb P}}

\title{$\SU_X(r,L)$ is separably unirational}
\date{\today}

\author[G. Hein]{Georg Hein}
\date{October 20, 2008}
\address{Universit\"at Duisburg-Essen, Fachbereich
Mathematik, 45117 Essen, Germany}

\email{georg.hein@uni-due.de}

\keywords{moduli space, vector bundles on a curve, separably
unirational}

\subjclass[2000]{14H60, 14D20}

\begin{abstract}
We show that the moduli space of $\SU_X(r,L)$ of rank $r$ bundles of
fixed determinant $L$ on a smooth projective curve $X$ is separably
unirational.
\end{abstract}

\maketitle

\section{Introduction}
In a discussion V.~B.~Mehta pointed out to me that for certain
applications about the cohomology of moduli of vector bundles on smooth
projective curves over algebraically closed fields in characteristic $p$
it is necessary to have that $\SU_X(r,L)$ is separably unirational.
This short note provides us with a proof of this statement.

\begin{localtheorem}\label{thm1}
Let $X$ be a projective curve of genus $g \geq 2$ over an algebraically
closed field $k$ of arbitrary characteristic. We fix a line bundle $L$
on $X$.  The moduli space $\SU_X(r,L)$ of $S$-equivalence classes of
semistable vector bundles of rank $r$ with determinant isomorphic to $L$
is separably unirational, that means there exists an open subset
$U \subset \pdop^{(r^2-1)(g-1)}$, and an \'etale morphism
$U \to \SU_X(r,L)$.
\end{localtheorem}
For a discussion of the notion {\em separable unirationality} and
typical applications see the lecture notes \cite[1.10 and 1.11]{Smi}.
If the characteristic of $k$ is zero, then separably unirational and
unirational coincide. Thus, in this case the result is well known
(see for example page 53 in Seshadri's lecture notes \cite{Ses}).

\section{Proof of theorem \ref{thm1}}
Let $d=\deg(L)$ be the degree of $L$. Fix a line bundle $M$ such that
$\deg(M) <\frac{d}{r}-2g$.
Let $E$ be any semistable vector bundle on $X$ with $\rk(E)=r$, and
$\det(E) \cong L$.
Set $M_0:=M^{\oplus (r+1)}$, and
$M_1:=M^{\otimes (r+1)} \otimes L\inv$.

If $\Ext^1(M,E) \ne 0$, then by Serre duality $\Hom(E, M \otimes
\omega_X) \ne 0$. $M \otimes \omega_X$ is a stable bundle
of slope $\mu(M \otimes \omega_X)= \deg(M)+2g-2 < \mu(E)$. Thus, there
can be only the zero morphism in $\Hom(E, M \otimes \omega_X) 0$.
So we have $\Ext^1(M,E) = 0$.

By the same argument we conclude that $\Ext^1(M,E(-P))=0$ for every
point $P \in X(k)$. Therefore $\Hom(M,E) \to \Hom(M, E\otimes k(P))$ is a
surjection.
We conclude that $\Hom(M,E) \otimes M \to E$ is surjective.
Since $X$ is of dimension one,
for a general subspace $V \subset \Hom(M,E)$ of dimension $r+1$ the
restriction $V \otimes M \to E$ is surjective. We obtain a surjection
$M_0 \rarpa{\pi} E$. The kernel of $\pi$ is a line bundle and its
determinant is $\det(\ker(\pi)) \cong M_0 \otimes \det(E)\inv =
M^{\otimes (r+1)} \otimes L\inv =M_1$. Taking one isomorphism
$M_1 \rarpa{\sim} \ker(\pi)$ we obtain for any $E$ as before the
existence of a short exact sequence
\begin{equation}\label{eq1}
0 \to M_1 \rarpa{\iota} M_0 \rarpa \pi E \to 0 \,.
\end{equation}
We use this to parameterize all bundles in $\SU_X(r,L)$ as cokernels of
morphisms $M_1 \to M_0$. To do so, we define $V:=\Hom(M_1,M_0)\dual$,
consider
$\xymatrix{\pdop(V) & \pdop(V) \times X \ar[l]_-p \ar[r]^-q & X}$,
and the natural morphism
$p^*\Ocal_{\pdop(V)}(-1) \otimes q^*M_1 \rarpa{ \alpha} q^*M_0$.  
We denote the cokernel of $\alpha$ by $\Ecal$.

Let $U_1$ be the open subset of points $u \in \pdop(V)$ such that
$\Ecal_u:= q_*(\Ecal \otimes p^*k(u))$ is a semistable bundle on $X$. We
see that the resulting morphism gives a surjection
$U_1 \rarpa \rho \SU_X(r,L)$.

Next we show that $\rho$ is infinitesimal
surjective. We take a point $[\iota] \in \pdop(V)$ corresponding to a
short exact sequence (\ref{eq1}). Let $D=k[\eps]/\eps^2$ be the ring of
dual numbers. To give an infinitesimal deformation of $E$ corresponds to
give a flat family $E_D$ on $X_D=\Spec(D) \times X$ which specializes to
$E$ when restricting to the reduced fiber $X_0 \cong X$.
Since we want to consider deformations with fixed determinant we have an
isomorphism $\det(E_D)\cong q_D^*L$ where $q_D$ is the projection
$X_D \to X$. The flat deformation $E_D$ yields a short
exact sequence $0 \to E \to E_D \to E \to 0$ on $\Spec(D) \times _k X$
which gives the exact sequence
\[ \xymatrix{\Hom_{X_D}(q_D^*M_0, E_D) \ar[r] &
\Hom_{X_D}(q_D^*M_0, E) \ar[r] \ar@{=}[d] & \Ext_{X_D}^1(q_D^*M_0, E)
\ar@{=}[d] \\ & \pi \in \Hom_X(M_0,E) & \Ext_X^1(M_0,E)} \]
Since $M_0=M^{\oplus(r+1)}$ we conclude from $\Ext^1(M,E)=0$ that
$\Ext^1(M_0,E)=0$. So $\pi$ is the
restriction of some $\pi_D \in \Hom_{X_D}(q_D^*M_0, E_D)$ to the reduced
fiber. The morphism $\pi_D:q_D^*M_0 \to E_D$ is surjective. Again
$\ker(\pi_D)$ is isomorphic to the line bundle $\det(q_D^*M_0) \otimes
\det(E_D)\inv = q_D^*(\det(M_0) \otimes L \inv) = q_D^*M_1$. Fixing such
an isomorphism we obtain a short exact sequence of $\Ocal_{X_D}$ bundles
\[ 0 \to q_D^* M_1 \rarpa{\iota_D} q_D^* M_0 \rarpa{\pi_D} E_D \to 0
\,.\]
We conclude that any deformation $E_D$ of $E$ is induced by a
deformation $\iota_D$ of $\iota$.
Now for a general linear subspace $L \subset \pdop(V)$ of dimension
$\dim\SU_X(r,L)=(r^2-1)(g-1)$ passing through a stable $[E] \in
\pdop(V)$ the composition of tangent maps is an isomorphism to
$T_{\SU_X(r,L),E}$. Thus, on some Zariski open subset $U \subset (L \cap
U_1)$ containing $[E]$ the morphism $U \to \SU_X(r,L)$ is \'etale.
\qed


\begin{thebibliography}{9}
\bibitem{Ses} C.~S.~Seshadri, {\em Fibr\'es vectoriels sur les courbes
alg\'ebriques}, Ast\'erisque {\bf 96}, Soci\'et\'e Math\'ematique de
France, Paris, 1982.
\bibitem{Smi} K.~Smith with an appendix by J.~Rosenberg,
{\em Rational and Non-Rational Algebraic
Varieties: Lectures of J\'anos Koll\'ar}, {\tt alg-geom/9707013}.
\end{thebibliography}
\end{document}